\documentclass[10pt,twocolumn,twoside]{IEEEtran}

\usepackage{srcltx,pdfsync}
\usepackage{amsmath}
\usepackage{amsfonts}
\usepackage{amssymb}
\usepackage{graphicx,subfigure,psfrag}
\usepackage{url,ifthen}
\usepackage[dvipsnames,usenames]{color}
\usepackage[colorlinks, linkcolor=Blue, filecolor =Red, citecolor=RoyalPurple,]{hyperref}

\usepackage[fancythm]{jphmacros2e}
\usepackage{etoolbox}
 
\newcommand{\G}{\mathcal{G}}
\newcommand{\V}{\mathcal{V}}
\newcommand{\N}{\mathcal{N}}
\newcommand{\A}{\mathcal{A}}
\newcommand{\B}{\mathcal{B}}
\newcommand{\M}{\mathcal{M}}
\newcommand{\MM}{\mathcal{S}_A}
\newcommand{\Mm}{\mathcal{S}_B}

\newcommand{\St}{\mathcal{S}}

\newcommand{\boldeps}{\mbox{$\boldsymbol{\epsilon}$}}

\newcommand{\boldp}{\mbox{$\boldsymbol{p}$}}

\newcommand{\boldx}{\mbox{$\boldsymbol{x}$}}
\newcommand{\boldone}{\mbox{$\boldsymbol{1}$}}
\newcommand{\boldy}{\mbox{$\boldsymbol{y}$}}
\newcommand{\boldz}{\mbox{$\boldsymbol{z}$}}
\newcommand{\boldu}{\mbox{$\boldsymbol{u}$}}

\newcommand{\R}{\mbox{$\mathbb{R}$}}

\newcommand{\barf}{\mbox{\underbar{$\nabla f$}}}

\def\ols{\ensuremath{\omega_\mathrm{L}}}%omega limit set

\newcommand{\jb}[1]{\footnote{{\bf JB's comment: #1}}}
\renewcommand{\jb}[1]{} %uncomment for final version

\begin{document}

%\begin{frontmatter}

\title{Distributed demand-side contingency-service provisioning while minimizing consumer disutility through local frequency measurements and inter-load communication} % Title, preferably not more 

%\thanks[footnoteinfo]{This paper was not presented at any IFAC 
%meeting. Corresponding author J.~Brooks. Tel. +1-407-617-5105. The research presented here was partially supported by the NSF through grant 1646229 and by the DOE BTO through a GMLC award titled ``Virtual Batteries''.}
                                            
\author{Jonathan~Brooks
        and Prabir~Barooah
\thanks{J. Brooks and P. Barooah are with the Department
of Mechanical and Aerospace Engineering, University of Florida, Gainesville,
FL, 32611 USA.}% <-this % stops a space
%\thanks{This paper was not presented at any IFAC meeting.}% <-this % stops a space
\thanks{Corresponding author J.~Brooks. e-mail: JonathanBrooksUF@gmail.com.}}
          
%\begin{keyword}                           % Five to ten keywords,  
%smart grid; distributed control; load control; contingency reserve.               % chosen from the IFAC 
%\end{keyword}                             % keyword list or with the 
                                          % help of the Automatica 
                                          % keyword wizard

\maketitle

\begin{abstract}                          % Abstract of not more than 200 words.
We consider the problem of smart and flexible loads providing contingency reserves to the electric grid and provide a Distributed Gradient Projection (DGP) algorithm to minimize loads' disutility while providing contingency services. Each load uses locally obtained grid-frequency measurements and inter-load communication to coordinate their actions, and the privacy of each load is preserved: only gradient information is exchanged---not disutility or consumption information. We provide a proof of convergence of the proposed DGP algorithm, and we compare its performance through simulations to that of a ``dual algorithm'' previously proposed in the literature that solved the dual optimization problem. The DGP algorithm solves the primal problem. Its main advantage over the dual algorithm is that it is applicable to convex---but not necessarily strictly convex---consumer disutility functions, such as a model of consumer behavior that is insensitive to small changes in consumption, while the dual algorithm is not. Simulations show the DGP algorithm aids in arresting grid-frequency deviations in response to contingency events and performs better or similarly to the dual algorithm in cases where the two can be compared. 
\end{abstract}

\IEEEpeerreviewmaketitle
%\end{frontmatter}

\section{Introduction}
\label{sec:intro}
For stable and reliable operation of the power grid, generation must match consumption at all timescales~\cite{kirby2007ancillary}.
Traditionally, this is achieved through controllable generators that adjust their outputs to track consumption. These services provided by generators are called ancillary services. %Contingency reserves are one such service provided after a sudden change in generation. 
% One form of ancillary services is contingency service, traditionally provided by fast ramping generation that can quickly change their generation in response to a contingency event, such as an unforeseen generation trip or transmission line trip. 
With the increasing penetration of volatile renewable energies into the power grid, more resources are required to provide these ancillary services. One form of ancillary services is contingency reserves, which are used to restore the consumption-generation balance after a sudden change (e.g., a large generator going off-line). % Greater need is felt in those regions of the world where renewable genration is seen to vary by hundreds of megawatts in only a few minutes~\cite{kamath2010understanding, is it appropriate?}. 
Conventional fossil-fuel generators are often operated at part-load in order to provide spinning reserves (fast-acting contingency reserves). However, rapidly ramping and operating at part-load can result in increased emission rates due to inefficient operation~\cite{lew2013western}. Part-loading requires additional generators to supply the needs of the grid as well. Building additional fossil-fuel generators to mitigate renewable volatility will reduce the environmental benefits of the renewable energies. 

It has been recognized in recent years that an attractive alternative exists: loads can be used to provide spinning reserves by changing their consumption without increasing emissions~\cite{milligan2010utilizing,siano2014demand}.
Due to the size of the grid, centralized solutions to the load control problem are not practical. A distributed solution is more attractive and is possible by utilizing the cyber-physical nature of the electric grid whereby ``information can be transmitted through actuation and sensing''~\cite{bolzam:13}. 
In particular, loads can provide primary control by using local frequency measurements~\cite{NERCBalancingDocument:2011short,pnnl2008value}. The value of information contained in frequency measurements has been recognized much earlier~\cite{schFAPER80}. Recent work in this vein includes~\cite{molina2011decentralized,short2007stabilization} in which loads are turned on or off based on frequency measurements.
Information from frequency measurements allows solutions not generally considered in the literature of distributed optimization (e.g.,~\cite{nedic2009approximate,zhu2012distributed,zhu2013approximate}).

Any changes in consumption to help the grid, however, may incur some cost or disutility for the consumer---such as deviation of the indoor temperature from a comfortable range. Thus, there is a need to balance the two---service to the grid and cost to the consumer. In this paper, we consider the problem of designing decision-making algorithms that provide spinning reserves through control of loads while striking this balance.
%This architecture, which is the focus of this paper, has the advantage of being fully distributed and scalable to a large number of loads. 

In this paper, we adopt the problem formulation from~\cite{zhao2013optimal}: minimize total consumer disutility while returning the consumption-generation mismatch in the grid to zero after a sudden change in generation. The consumption-generation mismatch is estimated by each load from noisy local frequency measurements using a state estimator. The algorithm in~\cite{zhao2013optimal} is based on solving the dual optimization problem, and we refer to it as the ``dual algorithm''. The dual variable, which is constant across the grid, is iteratively estimated using consensus averaging through inter-node communication.

\begin{figure}[ht]
  \centering
  \includegraphics[scale=0.325]{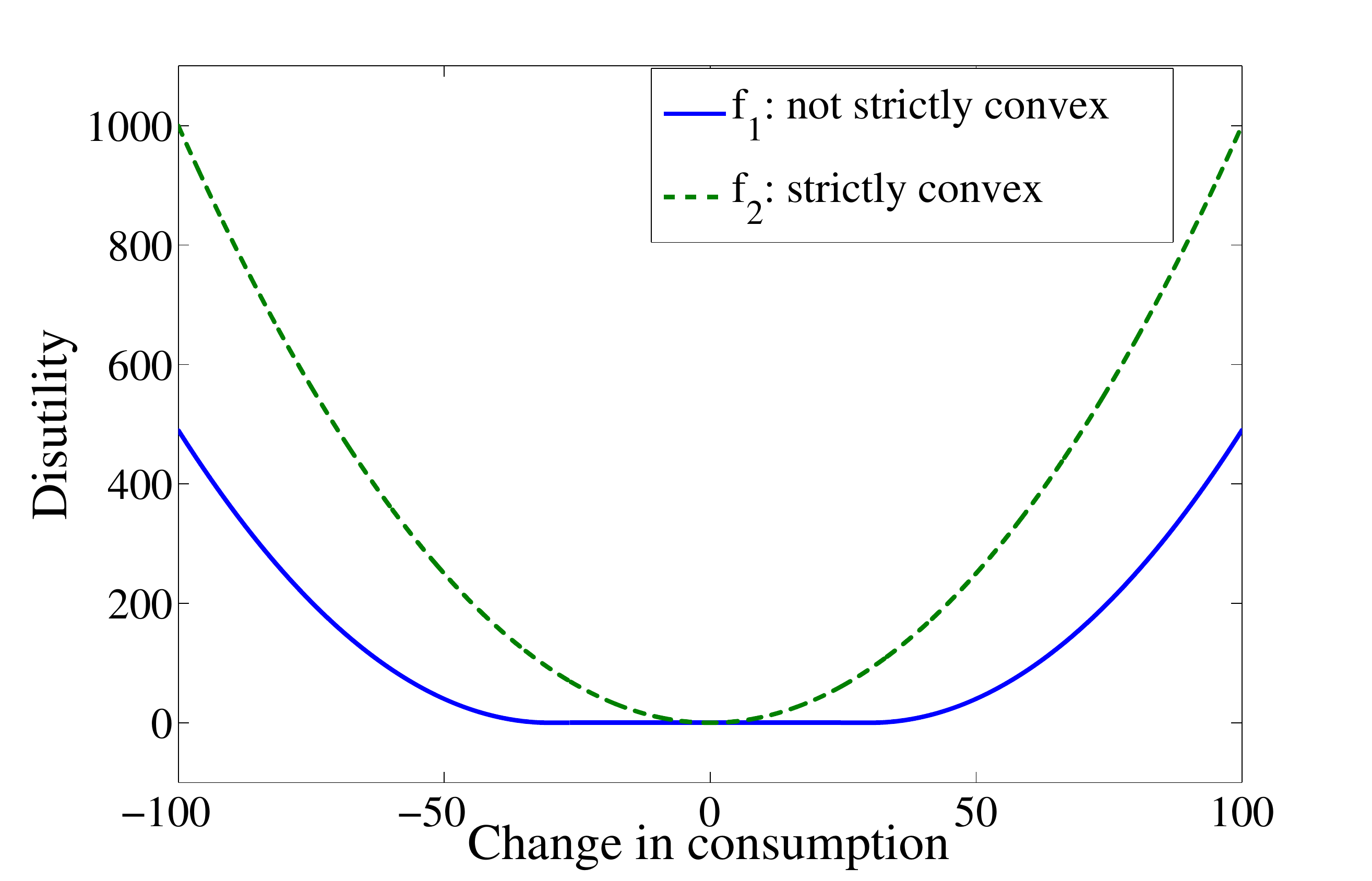}
  \caption{Two distinct models of consumer disutility.}
  \label{fig:disutility}
\end{figure}

The dual algorithm proposed in~\cite{zhao2013optimal} requires the consumers' disutilities to be strictly convex functions of changes in consumption. Quantifying consumers' disutility in response to consumption changes is challenging, and work in this area is limited. In~\cite{khadgi2014modeling}, an exponential function is used to model disutility, while~\cite{zhang2014optimal} proposes a dynamic disutility model. A study of an industrial aluminum-smelting plant suggests that there may be no disutility for several hours when changing consumption within some threshold of a nominal value, but there \emph{is} significant disutility if consumption is varied too much or for too long~\cite{todd2008providing}. Likewise,~\cite{linbarmeymid:2015} showed that consumption in commercial air-conditioning loads can be varied to provide ancillary services without any disutility (adverse effect on indoor climate) as long as the changes in consumption are small in amplitude and bandwidth-limited. Based on these studies, we hypothesize that an appropriate model of disutility for many consumers is like the function, $f_1$, shown in Figure~\ref{fig:disutility}. The disutility is zero for small changes in consumption but non-zero disutility for larger changes. Such a consumer's disutility is modeled by a convex---not strictly convex---function of consumption change.

This work is an extension of our prior work~\cite{brobar:ACC:2016}, where we proposed a method to solve the \emph{primal} problem in a distributed manner while preserving privacy, which we called the Distributed Gradient Projection (DGP) algorithm. The main advantage of the DGP algorithm over the dual algorithm of~\cite{zhao2013optimal} is that the DGP algorithm is applicable to disutility functions that are convex but not necessarily strictly convex.
The dual algorithm in~\cite{zhao2013optimal} required a strictly convex disutility function because the inverse image of the function's gradient is used in the computation, but the inverse image does not exist if there is a flat region in the disutility function (such as in $f_1$ in Figure~\ref{fig:disutility}). 

This paper makes several contributions over our preliminary work~\cite{brobar:ACC:2016}, in which the convergence proof was limited to the case where there are no upper and lower bounds on how much a load can change its consumption. That scenario is unrealistic and was adopted for tractability of analysis. In this work, we remove this assumption and prove convergence in the almost-sure sense as well as global asymptotic stability. The analysis of convergence in a bounded domain is substantially more involved because of the discontinuous nature of the resulting dynamics. The analysis uses the o.d.e. method of stochastic approximation~\cite{borkarbook:2008,kushner2003stochastic}. Connection to the so called Skorokhod problem helps resolve existence and uniqueness issues of the resulting o.d.e, but proving convergence to the desired set---solutions to the optimization problem---required novel analysis. We also present an example of when the algorithm does not converge---when the assumption of strict feasibility of the optimal solutions is violated, showing that the convergence result obtained here is the strongest one possible for the DGP algorithm. %In addition, simulation comparisons between the DGP and the dual algorithms presented here are more extensive compared to those in~\cite{brobar:ACC:2016}.

Simulations indicate that the DGP algorithm is effective in reducing frequency excursions following step changes in generation. Simulation comparisons, in those scenarios where comparison is possible, show that the proposed DGP algorithm performs better than or comparably to the dual algorithm of~\cite{zhao2013optimal}. On the other hand, our analysis is limited to time-invariant communication, whereas convergence for the dual algorithm was proved to hold even for time-varying communication~\cite{zhao2013optimal}.

This paper is organized as follows. Section~\ref{sec:statement} formally defines the problem that we solve. In Section~\ref{sec:sgp}, we describe the DGP algorithm. We present the convergence result in Section~\ref{sec:convergence} and proofs in Section~\ref{sec:proofs}. Simulations are described in  Section~\ref{sec:results}. Finally, Section~\ref{sec:conc} concludes this work and discusses avenues for future extensions.

\section{Problem Formulation}
\label{sec:statement}
As in~\cite{zhao2013optimal}, we consider an electric grid with a single frequency throughout the grid, whose nominal value is denoted by $\omega^*$, such as in a microgrid. %This is the case when electrical distances are negligible---such as in a microgrid. 
There are $n$ controllable loads. The deviation of load $i$'s consumption from its nominal value is denoted by $x_i$ and incurs a disutility, $f_i(x_i)$. The deviation must lie in $\Omega_i \triangleq [$\b{$x$}$_i,\bar{x}_i]$, specified a-priori. 

Let $\Delta g$ be the generation deviation from the nominal value. The problem is for the loads to decide how much to change their own consumption so that the consumption-generation mismatch is diminished while the resulting total disutility of the loads is minimized:
\begin{align}\label{eq:primal}
\min_{x_i,\; i=1,\dots,n}\sum_{i=1}^{n}f_i(x_i), \;\; \mathrm{s.~t.~}\sum_{i=1}^{n}x_i = \Delta g, \quad x_i \in \Omega_i,
\end{align}

Load $i$ can obtain a noisy measurement of the grid frequency and can use it to make a decision on $x_i$. In addition, the computation of the decision variables, $x_i$, must be distributed in the following sense. There is a communication graph, $\G=(\V,\mathcal{E})$, where the node set, $\V=\{1,2,\dots,n\}$, is simply the loads and the edge set, $\mathcal{E} \subset \V \times \V$, specified a-priori, determines which pairs of loads can exchange information. The set of neighbors of load $i$, with which it can exchange information, is defined by $\N_i \triangleq \{\;j\;|\;(i,j)\in\mathcal{E}\}$, and we denote $n_i\triangleq|\N_i|$ as the number of neighbors of load $i$. 
The frequency measurements are essential since every load can use them to estimate the equality-constraint violation, $u \triangleq \Delta g - \sum_{i=1}^{n}x_i$. How this is done is described in Section~\ref{sec:estimator}.

%Although Problem~\eqref{eq:primal} does not include time, time plays a role since the noise on frequency measurement is naturally modeled as a stochastic process, and consequently the estimates of $u$ obtained by every node vary with time.\pb{also, power system dynamics will cause transients. iterative optimization will use estimate of mismatch which will vary with iteration since new measurements that vary due to transients and noise will affect it. Simplify to emphasize the most important ones.}
Although Problem~\eqref{eq:primal} does not include time, the proposed algorithm is an iterative approach, so time does play a role. The reason for this is twofold: i) due to the limited information available to the loads, they are unable to find the solution in just one step, so they move ``toward'' the solution and then re-evaluate at the new location ii) due to generator dynamics, the frequency, $\omega$, is dependent on the loads, so as the loads change in time, so does the frequency, which the loads can then use for feedback control.
Time is measured by a discrete iteration counter: $k=0,1,\dots$. The generation at time $k$ is denoted by $g[k]$ so that the generation change from nominal is $\Delta g[k] \triangleq g[k] - g^*$, where $g^*$ is the nominal generation. We  assume that, at $k=0$, total load and total generation are equal, and we limit ourselves to step changes in generation. That is, $\Delta g[k]= 0$ for $k<K$ for some $K$, and $\Delta g[k] = \bar{g}$ for $k\geq K$, where $\bar{g}$ is the size of the step change.

\section{Distributed Gradient Projection (DGP) Algorithm}
\label{sec:sgp}
To describe the algorithm, we define the consumption-generation mismatch at iteration $k$:
\begin{align}\label{eq:u}
 u[k] \triangleq \Delta g[k] - \sum_{i=1}^{n}x_i[k] = \Delta g[k] - \boldone^T\boldx[k], 
\end{align}
where $\boldx[k] \triangleq [x_1[k],\dots,x_n[]^T $ and $\boldone \in \R^n$ is a vector of all ones. The frequency deviation from its nominal value is denoted by $\Delta\omega[k]\triangleq\omega[k]-\omega^*$. Neither $\sum_{i=1}^{n} x_i[k]$ nor $\Delta g[k]$ is known to any of the loads. However, 
load $i$ can obtain a noisy measurement, $\Delta\tilde{\omega}_i[k]$, of the frequency deviation. It uses this measurement to estimate the mismatch, which is denoted by $\hat{u}_i[k]$, 

The update law of the DGP algorithm comprises of 3 main operations: (i) a generation-matching step, (ii) a gradient descent step, and (iii) a projection step. The first step uses the estimated mismatch, $\hat{u}_i[k]$, to compute a consumption change that will reduce the mismatch. For the second step, pure gradient descent, though possible, will violate the equality constraint (consumption-generation matching). Therefore the gradient descent step is designed to be orthogonal to the generation-matching step; i.e., it does not change the total consumption. The updates computed by the first two steps are added, and their sum is then projected onto $\Omega_i$ to respect the upper and lower bounds on consumption change.

The update law of the DGP algorithm at load $i$ at time $k$ is summarized below:

\textbf{\emph{DGP Algorithm}}:
\begin{enumerate}
\item Obtain an estimate of $u_i[k]$ (call it $\hat{u}_i[k]$) from the measurement, $\Delta\tilde{\omega}_i[k]$, using a state estimator (described in Section~\ref{sec:estimator}). 
\item Compute gradient $\frac{d}{dx_i} f_i(x_i[k])$, transmit gradient value to neighbors, receive neighbors' gradient values, and compute 
\begin{align}\label{eq:deltax-0}
\Delta x_i[k] := - n_i\nabla f_i(x_i[k]) + \sum\limits_{j\in\N_i}\nabla f_j(x_j[k]).
\end{align}
\item Update load change as $x_i[k+1] = P_{\Omega_i}\big[x_i[k] + \alpha[k]\Delta x_i[k] + \gamma[k]\hat{u}_i[k]\big]$, where $P_{\Omega_i}[\cdot]$ denotes the standard projection operator and $\alpha[k],\gamma[k]$'s are step sizes.\break
\end{enumerate}

\noindent Note that loads only exchange gradient information, so each load's disutility and consumption information remain private and are not shared with other loads.

\subsection{Estimation of consumption-generation mismatch using frequency measurements}\label{sec:estimator}
We borrow the estimation method proposed in~\cite{zhao2013optimal} for estimating $u$, although it is possible to use any estimator in the DGP algorithm. The power grid is modeled as a discrete-time LTI system with consumption-generation mismatch, $u[k]$, as the input and frequency deviation from nominal, $\Delta\omega[k]$, as the output. At each time $k$, load $i$ obtains the noisy measurement, $\Delta\tilde{\omega}_i[k]$, to estimate the state of the plant by using the estimator in~\cite{kitanidis1987unbiased}, which was developed for estimating the state of a system with an unknown input. 
Once the state estimate is obtained, each load estimates the unknown input by effectively assuming that the most recent output is error-free and then solving for the previous input from the state equations. 

We denote the estimation error at time $k$ by 
\begin{align}\label{eq:epsilon-def}
  \boldeps[k]\triangleq\hat{\boldu}[k]-u[k]\boldone,
\end{align}
where $\hat{\boldu}[k]$ is the column vector of $\hat{u}_i[k]$'s. Define the $\sigma$-algebra, $\mathcal{F}[k-1]:=\sigma(\epsilon_i[\ell-1]~|~i\in\V;\;1\leq \ell \leq k)$. It was shown in~\cite{kitanidis1987unbiased} that
\begin{align}\label{eq:error}
\mathbb{E}\big[\epsilon_i[k]|\mathcal{F}[k-1]\big]=0.
\end{align}

The following proposition is reproduced from~\cite{zhao2013optimal}.
\begin{proposition}[\cite{zhao2013optimal}]\label{prop:variance}
Let $A$, $B$, and $C$ denote the process, input, and output matrices of the state-space model of the power grid used in the state estimator described in~\cite{kitanidis1987unbiased}, and let $I_A$ be an identity matrix the size of $A$. If every eigenvalue of $(I_A-B(CB)^{-1}C)A$ lies within the unit circle, then $\lim\limits_{k\to\infty}\mathbb{E}[(\epsilon_i[k])^2|\mathcal{F}[k-1]\big]$ exists.
\end{proposition}
The following corollary is a straightforward consequence of~\eqref{eq:error} and Proposition~\ref{prop:variance} based on the definition of a martingale-difference sequence. 
\begin{corollary}\label{cor:martin}
If the condition for Proposition~\ref{prop:variance} holds, then~\eqref{eq:error} and Proposition~\ref{prop:variance} imply that the estimation error sequence, $\boldeps[k]$, is a martingale-difference sequence. %, which satisfies Assumption~\ref{as:tech}(3) below.
\end{corollary}

\section{Convergence analysis}\label{sec:convergence}
We make the following assumptions for our analysis.

\begin{assumption}\label{as:f}
\textbf{\textup{(Assumptions on disutility).}}
\begin{enumerate}
\item $f_i(x_i)$ is convex for each $i$ with a (not necessarily unique) minimum at $x_i=0$.
\item $f_i(x_i)$ is coercive for each $i$; i.e., $\{x_i|f_i(x_i)\leq F\}$ is compact for every $F\geq 0$ for each $i$.
\item $f_i(x_i)$ is continuously differentiable for each $i$.
\item $\nabla f_i(x_i)$ is Lipschitz for each $i$.
\end{enumerate}
\end{assumption}

\begin{assumption}\label{as:power}
\textbf{\textup{(Geometric assumptions).}}
\begin{enumerate}
\item The domain, $\Omega$, is compact.
\item The communication graph, $\G$, is connected.
\item The disturbance is a constant: $\Delta g[k] \equiv \bar{g}$ for all $k\geq 0$.
\end{enumerate}
\end{assumption}

\begin{assumption}\label{as:tech}
\textbf{\textup{(Technical assumptions).}}
  \begin{enumerate}
  \item $\alpha[k] = c\gamma[k]$ for some positive constant $c$.
  \item The function, $\gamma[k]\to 0$, satisfies $\sum_{k=0}^{\infty}\gamma[k]=\infty$ and $\sum_{k=0}^{\infty}(\gamma[k])^2< \infty$.
  \item The estimation error sequence, $\boldeps[k]$, is a martingale-difference sequence.
  \end{enumerate}
\end{assumption}

Assumption~\ref{as:f} is readily satisfied because the disutility functions are a modeling choice. Assumption~\ref{as:power}(1) is always met in practice since a load cannot change its demand outside its maximum rated power and 0. Assumption~\ref{as:power}(3) means that we only consider a step-change in generation, which is a good approximation of a contingency event. Assumptions~\ref{as:tech}(1-2) are satisfied by choosing $\alpha[k]$ and $\gamma[k]$ appropriately. Assumptions~\ref{as:tech}(2-3) are standard in the field of stochastic approximation~\cite{borkarbook:2008}.
Note that Assumption~\ref{as:tech}(3) holds in this work due to Proposition~\ref{prop:variance}.
%Assumption~\ref{as:tech}(2-3) are standard in the field of stochastic approximation~\cite{borkarbook:2008}.  Assumption~\ref{as:tech}(1-2) can be satisfied by choosing $\alpha[k]$ and $\gamma[k]$ appropriately, while the third holds in this work due to Proposition~\ref{prop:variance}. %Assumption~\ref{as:f}(1) is always met in practice since a load cannot change its demand outside its maximum rated power and 0. 
%The properties of $f(\cdot)$ in Assumption~\ref{as:f} are readily satisfied because the disutility functions are a modeling choice. Assumption~\ref{as:power} means that we only consider a step-change in generation, which is a good approximation of a contingency event.

The main convergence result is the following. 
\begin{theorem}\label{the:main}
If Assumptions~\ref{as:f},~\ref{as:power}, and~\ref{as:tech} hold and all solutions to Problem~\eqref{eq:primal} are strictly feasible, $\boldx[k]$ converges to a solution to Problem~\eqref{eq:primal}, almost surely (a.s.).
\end{theorem}

The proof of this result relies on the so called o.d.e. method of stochastic approximation, which establishes a rigorous connection between noisy discrete iterations and a continuous-time o.d.e.~\cite{borkarbook:2008,kushner2003stochastic}, stated next.
\begin{proposition}[Theorem 2.1 (Chapter 5) in~\cite{kushner2003stochastic}]\label{prop:borkar}
Consider the sequence $\{\boldy[k]\}$ generated by the iteration
\begin{align*}
\boldy[k+1] = P_\Omega\Big[\boldy[k] + \gamma[k]\big(h(\boldy[k])+\boldeps[k]\big)\Big],
\end{align*}
where $P_\Omega$ is the projection operator onto $\Omega$, $h(\boldy):\R^n\to\R^n$ is Lipschitz, and $\{\boldsymbol{\epsilon}[k]\}$ is a martingale-difference sequence. If $\gamma[k]$ satisfies Assumption~\ref{as:tech}(2) then $\boldy[k]$ converges a.s. to some limit set of the o.d.e.,
\begin{align*}
\dot{\boldy}(t)=\Gamma_{\Omega,\boldy(t)}[h(\boldy(t))],
\end{align*}
where $\Gamma_{\Omega,\boldy(t)}[\cdot]$ denotes the continuous-time projection operator of $\boldy(t)\in\R^n$ onto $\Omega$. That is, the $i^\mathrm{th}$ component of $\Gamma_{\Omega,\boldy}[\boldz]$ is
\begin{align*}
\big(\Gamma_{\Omega,\boldy}[\boldz]\big)_i = \left\{\begin{array}{l} 0;~~y_i=\min \Omega_i,~z_i<0\\ 0;~~y_i=\max \Omega_i,~z_i>0\\ z_i;~~\mathrm{o.w.} \end{array} \right.
\end{align*}
\end{proposition} 
By using~\eqref{eq:epsilon-def} and Assumption~\ref{as:tech}, the update law of the DGP algorithm can be written as 
\begin{align}\label{eq:i-DGP-connect2ode}
\begin{aligned}
x_i[k+1] &= P_{\Omega_i}\bigg[x_i[k]\\
         &\qquad\quad + \gamma[k]\Big(c\sum_{j\in\N_i} \big(\nabla f_j(x_j[k]) - \nabla f_i(x_i[k])\big)\\
         &\qquad\qquad + u[k] + \epsilon_i[k]\Big)\bigg],
\end{aligned}
\end{align}
or compactly,
\begin{align}\label{eq:DGP-connect2ode}
\begin{aligned}
\boldx[k+1] = P_{\Omega}&\big[\boldx[k]\\
&+ \gamma[k](-cL\nabla f(\boldx[k])^T+ u[k]\boldone+\boldeps[k])\big],
\end{aligned}
\end{align}
where $\nabla f(\boldx)$ is the gradient of $f(\boldx):=\sum_{i=1}^n f_i(x_i)$ and $L \in \R^{n \times n}$ is the graph Laplacian of $\G$~\cite{GodsilRoyle_2001}: $L_{ii}=n_i$ (number of neighbors of node $i$), and for $j\neq i$, $L_{ij}=-1$ if $j \in \N_i$, $L_{ij}=0$ if $j \notin \N_i$. By Proposition~\ref{prop:borkar}, the iterates, $\boldx[k]$, converge a.s. to a limit set of the o.d.e.,
\begin{align}
\label{eq:ode}
\dot{\boldx}(t) &=\Gamma_{\Omega,\boldx(t)}[-cL\nabla f(\boldx(t))^T + u(t)\boldone],
\end{align}
\begin{align}
\label{eq:u(t)-def}
u(t) &:=\bar{g}-\boldone^T\boldx(t)
\end{align}
The task is to prove that the trajectories of o.d.e.~\eqref{eq:ode} converge to the set of solutions to  Problem~\eqref{eq:primal}, which we denote by $X^*$. Because of the convexity, $\boldx^* \in X^*$ if and only if it satisfies the first-order necessary conditions of optimality~\cite{Luenbergerbluebook:2003}. We denote the boundary of $\Omega$ by $\partial\Omega$, and we denote the interior of $\Omega$ by $\Omega^o$. Because of the hypothesis about strict feasibility in Theorem~\ref{the:main}, only the equality constraint is active, which is also regular, and thus it is straightforward to verify that 
\begin{align}\label{eq:Xstar}
  X^*=\{\boldx \in \Omega^o~|~\boldone^T\boldx=\bar{g},~\nabla
  f_i(x_i)=\nabla f_j(x_j),~\forall i,j\in\V\}.
\end{align}
Because each $f_i$ is convex, each gradient $\nabla f_i(x_i)$ is nondecreasing, and the optimal set, $X^*$, is connected. It follows that the gradient, $\nabla f(\boldx)$, at any optimal point $\boldx^*\in X^*$ is unique, which can be seen via contradiction. Suppose $\boldx^*_1,\boldx^*_2\in X^*$ with $\|\nabla f(\boldx^*_1)\|<\|\nabla f(\boldx^*_2)\|$%; recall that all elements of $\nabla f(\boldx^*_1)$ are equal by optimality, and the same applies for $\nabla f(\boldx^*_2)$
. Then $f(\boldx^*_1)<f(\boldx^*_2)$ by convexity, so $\boldx^*_2$ is not optimal---a contradiction. Therefore, $\nabla f(\boldx^*)$ is unique for all $\boldx^*\in X^*$. We denote the corresponding scalar by $\nabla f^*$:
\begin{align}\label{eq:gradfstar}
  \nabla f^* \triangleq \nabla f_i(x_i), \text{ where } \boldx \in X^* \text{ and } i \in \{1,\dots,n\}.
\end{align}
The constant, $\nabla f^*$, is called the \emph{optimal gradient} in the sequel.
% Note that $\nabla f^*$ is the same no matter which optimal $\boldx$ or which component $j$ is used\pb{I introduced this definition when I realized this quantity is used in a proof without defining it. Check if what I wrote is correct. But my def also feels too cumbersome. Compress? \textbf{JB:} We define it in the proof of Lemma 4.7. Perhaps we should just reproduce that definition here? It is essentially ``by convexity, the optimal gradient is unique. Denote that gradient by $\nabla f^*$.}. 

Since all limit sets of o.d.e~\eqref{eq:ode} are contained within its $\omega$-limit set, which we call $\ols$, it suffices to show that the optimal set, $X^*$, contains the $\omega$-limit set: $\ols\subseteq X^*$. Theorem~\ref{the:main} then follows immediately by Proposition~\ref{prop:borkar}.

In many applications of the o.d.e. method, the main hurdle in analyzing convergence is to establish boundedness of the iterates, $\boldx[k]$~\cite{borkarbook:2008}. In our case, presence of the projection step guarantees boundedness trivially. However, this projection makes the limit sets of the o.d.e. difficult to characterize due to the discontinuous nature of projected dynamical systems~\cite{nagurney2012projected}. In fact, existence and uniqueness of a solution to~\eqref{eq:ode} needs to be established first. Because the right-hand side of~\eqref{eq:ode} is discontinuous, existence of classical solutions to the o.d.e. (i.e., solutions that are continuously differentiable) is not guaranteed. For discontinuous dynamics, the notions of Caratheodory and Filippov solutions are used. Caratheodory solutions are essentially a generalization of classical solutions that are absolutely continuous and differentiable almost everywhere, while Filippov solutions are absolutely continuous maps that satisfy a differential inclusion almost everywhere~\cite{cor:2008}. When dealing with Filippov solutions, the way in which the differential inclusion is defined plays a crucial role in both existence and uniqueness. This can further complicate the analysis.

The o.d.e.~\eqref{eq:ode} is closely related to the so called Skorokhod problem, which deals with stochastic differential equations with boundaries~\cite{skorokhod1961stochastic}. The connection with the Skorokhod problem and ours is useful; known results on the Skorokhod problem provide guarantees of not only existence but also uniqueness of Caratheodory solutions to o.d.e.~\eqref{eq:ode}~\cite{nagurney2012projected}. For completeness, the relevant result from~\cite{nagurney2012projected} is presented below.

\begin{proposition}[Theorem 2.5 in~\cite{nagurney2012projected}]\label{prop:nag}
Let $h(\cdot)$ be Lipschitz continuous and
\begin{align}\label{eq:proj_ode}
\dot{\boldx}(t)=\Gamma_{\Omega,\boldx(t)}\big[h(\boldx(t))\big],
\end{align}
where $\Gamma$ is the projection defined in Proposition \ref{prop:borkar} and $\Omega$ is compact. For any $\boldx(0)\in\Omega$, there exists a unique Caratheodory solution, $\boldx(t)$, to~\eqref{eq:proj_ode} starting from $\boldx(0)$. Furthermore, $\boldx(t)$ is continuous with respect to the initial condition, $\boldx(0)$.
\end{proposition}

The second challenge is analyzing limiting behavior of the o.d.e. The trajectory may ``evolve[] along a 'section' of [$\partial\Omega$]....  At a later time the solution may re-enter [$\Omega^0$], or it may enter a lower[-]dimensional part of [$\partial\Omega$].''~\cite{nagurney2012projected}. The trajectory may go on doing so without ever converging to a limit set with tractable structure. This is a severe hurdle in analyzing dynamical systems evolving in a bounded region, where the boundedness is enforced through a projection operator. In~\cite{nagurney2012projected}, this hurdle was sidestepped by \emph{assuming} that the $\omega$-limit set is contained in the set of fixed points. That is equivalent to assuming that complicated limit sets such as limit cycles do not arise. The reference~\cite{borkarbook:2008} also mentions that, in presence of the projection operator onto domains with non-smooth boundaries (such as in our case), existence of an o.d.e. limit of the discrete-time algorithm may be a non-trivial issue. The classic reference~\cite{kushner2003stochastic} on stochastic approximation and o.d.e. methods mentions that basic results on convergence to a limit set (such as Proposition~\ref{prop:borkar}) is rarely useful since the limit set can be the whole of the bounded domain in which the o.d.e. evolves. This was also the main argument in \cite{nagurney2012projected} that is mentioned above. Only in the special case when the right-hand side of the o.d.e. is a gradient-descent system can convergence to a set of fixed points be established~\cite{kushner2003stochastic}. However,~\eqref{eq:ode} is not a gradient-descent system.  

Our main convergence result, Theorem \ref{the:main}, is obtained without making the strong assumptions that are typical, such as assuming that no limit cycles exist. The structure of the  Laplacian matrix plays an important role in the proofs of the main theorem and the technical results needed for the theorem proof. A key technical result is that the equality constraint is eventually satisfied, which is stated in the next lemma. 

\begin{lemma}\label{lem:u}
If Assumptions~\ref{as:f},~\ref{as:power}, and~\ref{as:tech} hold and all solutions to Problem~\eqref{eq:primal} are strictly feasible, $u(t)\to0$, where $u(t)$ is defined in \eqref{eq:u(t)-def} and $\boldx(t)$---in the definition of $u(t)$---is governed by o.d.e.~\eqref{eq:ode}.
\end{lemma}
A similar result was obtained in~\cite{cherukuri2014initialization}. However, the analysis in that work uses a penalty method in place of projection. Consequently, the result of Lemma~\ref{lem:u} is obtained immediately. In our case, where we retain the projection operation, the proof of Lemma~\ref{lem:u} is not as straightforward. We provide our proof of the lemma in Section~\ref{sec:proofs}. We are now ready to prove Theorem~\ref{the:main}. 

\noindent\textbf{Proof of Theorem~\ref{the:main}:}
Recall the discussion after Proposition~\ref{prop:borkar}: we must prove that the trajectories of o.d.e.~\eqref{eq:ode} converge to the set, $X^*$. By Proposition~\ref{prop:nag}, the solution to o.d.e.~\eqref{eq:ode} is continuous with respect to the initial condition, so $\ols$ is an invariant set~\cite{Khalil}. Our goal is to use this invariance property to characterize $\ols$ and show that $\ols \subseteq X^*$.
For convenience, we write the $i^\mathrm{th}$ element of~\eqref{eq:ode} below:
\begin{align}
\label{eq:iode}
\begin{aligned}
\dot{x}_i(t) = \Gamma_{\Omega_i,x_i(t)}\Big[& c\sum_{j\in\N_i} \Big(\nabla f_j(x_j(t)) - \nabla f_i(x_i(t))\Big) \\
& + u(t)\Big].
\end{aligned}
\end{align}
By Lemma~\ref{lem:u}, $\boldx(t)\to \{\boldx\in\Omega~|~\boldone^T\boldx=\bar{g}\}=:U$; this implies that all limit sets are contained within $U$. Because the $\omega$-limit set is the union of all limit points, we have $\ols \subseteq U$. Let 
\begin{align}\label{eq:gradf_lowbar}
  \barf(t)\triangleq\min_i\nabla f_i(x_i(t)).
\end{align}
Imagine a trajectory starting in $\ols$: $\boldx(0)\in\ols \subseteq U$, which implies $\boldx(t)\in U$ for all $t$ by invariance of $\ols$. Since $u(t)\equiv 0$ for this trajectory, o.d.e.~\eqref{eq:ode} reduces to $\dot{\boldx}(t) = \Gamma_{\Omega,\boldx(t)}\big[-cL\nabla f(\boldx(t))^T\big]$. That is, for every $i\in\V$,  $\dot{x}_i(t)$ is the average of its neighbors' gradients at $t$, which shows that $\dot{x}_{r(t)}\geq0$, where $r(t):=\arg\min_i\nabla f_i(x_i(t))$. Therefore, $\barf(t)=\nabla f_{r(t)}(x_{r(t)}(t))$ is nondecreasing by convexity. (Note that this does not require uniqueness of $r(t)$, so when two gradients are both minimal, both are nondecreasing.) Because $\barf(t)$ is bounded as well, which comes from boundedness of the domain, $\Omega$, $\barf(t)$ converges. Denote the limit by $F$ so that trajectories starting in $\ols$ converge to $\mathcal{F}:=\{\boldx\in\Omega~|~\min_i\nabla f_i(x_i)=F;~\boldone^T\boldx=\bar{g}\}\subseteq \ols \subseteq U$. Therefore, for any $\boldx(0)\in\ols$, $\boldx(t)\to\mathcal{F}$.

We now show that $\boldx(t)\to\mathcal{F}$ for any $\boldx(0)\in\Omega$. Suppose we have some trajectory, $\boldx(t)$ with $\boldx(0)=\boldx_o\in\Omega$. From the definition of the $\omega$-limit set, there exists a point, $\boldy\in\ols$, and a sequence, $t_0,t_1,\ldots$ with $t_n \to \infty$ as $n \to \infty$, such that $\boldx(t_n)\to\boldy$ as $n\to\infty$. Let $\boldx_n(t) := \boldx(t_n + t)$ for every $t$ so that $\boldx_n(0)=\boldx(t_n)$. Let $\boldx_{\ols}$ be a trajectory that starts from $\boldy$: $\boldx_{\ols}(0)=\boldy\in\ols$. By construction $\lim_{n\to\infty}\boldx_n(0)=\boldy = \boldx_{\ols}(0)$. Therefore, for any $T\in[0,\infty)$, $\boldx_n(T)\to\boldx_{\ols}(T)$ as $n\to\infty$ by continuity of $\boldx(t)$ with respect to $\boldx(0)$ (see~Proposition~\ref{prop:nag}). Now, we have
$\lim_{n\to\infty}\boldx(t_n+T)=\lim_{n\to\infty}\boldx_n(T)=\boldx_{\ols}(T)$. This is true for any $T\geq0$, so we may take the limit as $T\to\infty$ to obtain $\boldx(t)\to\mathcal{F}$ because $\boldx_{\ols}(T)\to\mathcal{F}$ as $T\to\infty$. Therefore, $\boldx(t)\to\mathcal{F}$ for any $\boldx(0)\in\Omega$, so $\ols\subseteq\mathcal{F}$. Since $\mathcal{F}\subseteq\ols$ also, we have $\mathcal{F}=\ols$. Therefore, $\mathcal{F}$ is invariant by invariance of $\ols$.

Next, we show that $\mathcal{F}$ is invariant only if $F=\nabla f^*$. Consider again a trajectory starting in $\mathcal{F}$: $\boldx(t_0)\in \mathcal{F}$ for some $t_0$. If $\nabla f_i(x_i(t_0))=F$ for every $i\in\V$, then the first-order necessary conditions for optimality are satisfied, so $F=\nabla f^*$, and $\mathcal{F}=X^*$, and the proof is complete. Suppose there exists an $\ell\in\V$ such that $\nabla f_\ell(x_\ell(t_0))>F$. 
Let $q:=\arg\min_i\nabla f_i(x_i(t_0))$; then $\dot{x}_q(t_0)\geq0$ with equality if and only if $\nabla f_j(x_j(t_0))=F$ for all $j\in\N_q$, in which case $q$ is not unique. Denote the different choices for $q$ by $q_i,~i\in\V$. 
Because $\nabla f_\ell(x_\ell(t_0))>F$ and the graph $\G$ is connected, $\dot{x}_{q_j}(t_0)>\epsilon$ for some $q_j$ and some $\epsilon>0$. It follows that there exists some $t_1>t_0$ such that $\nabla f_{q_j}(x_{q_j}(t_1))>F$. Note that $u(t)=0$ implies $\nabla f_{q_j}(x_{q_j}(t))>F$  for all $t\geq t_1$; this is because $\dot{x}_{q_j}(t)\to\R_{\geq0}$ as $\nabla f_{q_j}(x_{q_j}(t))\to F$ (see~\eqref{eq:iode}). Hence, $\dot{x}_{q_i}(t_1)>\delta$ for some $q_i\in\N_{q_j}$ and some $\delta>0$. By continuing this argument, we see that eventually there exists some $t_2>t_0$ such that $\barf(t_2)>F$; that is, when the gradient of one $q_j$ increases, that causes the gradient of all $q_i\in\N_{q_j}$ to eventually increase, and so on. However, $\barf(t_2)>F$ implies $\boldx(t_2)\notin\mathcal{F}$, so $\mathcal{F}$ is not invariant---a contradiction. Therefore, if $\boldx(t_0)\in\mathcal{F}$, there exists no $\ell\in\V$ such that $\nabla f_\ell(x_\ell(t_0))>F$. Hence, $\nabla f_i(x_i(t_0))=F$ for all $i\in\V$, which satisfies the optimality conditions, so $F=\nabla f^*$, and $\ols=\mathcal{F}=X^*$%. Thus, $F=\nabla f^*$, and $\mathcal{F}=X^*$. Therefore, $\ols\subseteq X^*$
, which completes the proof.\hfill$\blacksquare$

The following lemma shows that the optimal solution set, $X^*$, is not only globally attractive, but also globally asymptotically stable. The notion of stability of an invariant set used here is the one from~\cite{Khalil}. %The fact that $X^*$ is an invariant set of o.d.e.~\eqref{eq:ode} is trivial to establish.

\begin{lemma}\label{lem:stable}
Let the conditions of Theorem~\ref{the:main} hold. Then the optimal set, $X^*$, is globally asymptotically stable.
\end{lemma}

The proof of this lemma is provided in Section~\ref{sec:proofs}.\\

\noindent\textbf{Comment on strict feasibility:}
The assumption of strict feasibility of all solutions to Problem~\ref{eq:primal} is necessary for the convergence of the DGP algorithm to hold. If $X^*$ is not strictly feasible, then the DGP algorithm is not guaranteed to converge to $X^*$; in fact it may converge to non-optimal points. This may be seen through the following 2D counterexample. Let $\boldx\in\R^2$, $\bar{g}=1$, $f_1(x_1)=(x_1)^2$, $f_2(x_2)=(x_2)^2$, $\Omega = [0, \;  1/4] \times [0,\; 1]$, and $c=1$. $\G$ consists of the two nodes and one edge connecting them. The constraint, $u=0$, is satisfied on the line $x_1+x_2=1$. It is straightforward to verify that $\boldx^*=[1/4,~3/4]^T$ is the unique solution to Problem~\eqref{eq:primal}, but $\boldx^*$ is not strictly feasible since it lies on the boundary of $\Omega$. Now, from \eqref{eq:ode},
\begin{align*}
\dot{x}_2(t)|_{\boldx^*}&=\Gamma_{\Omega_2,x_2(t)}[-2x_2(t)+2x_1(t)+u(t)]|_{\boldx^*}\\
            &=\Gamma_{\Omega_2,0.75}[-2\cdot0.75 + 2\cdot0.25 + 0]\\
            &=-1,
\end{align*}
which shows that $[1/4,~3/4]^T$ is not an equilibrium point of o.d.e. \eqref{eq:ode}. Therefore, the DGP algorithm is not guaranteed to converge to $X^*$. It can also be shown that for this example $[1/4,~5/12]^T$ is an equilibrium point of o.d.e. \eqref{eq:ode} and that it is attractive, but $[1/4,~5/12]^T$ is not a solution to Problem \eqref{eq:primal}. Therefore, from Proposition \ref{prop:borkar}, the iterates will converge almost surely to a non-optimal point.

\subsection{Proofs}\label{sec:proofs}
\noindent \textbf{Proof of Lemma~\ref{lem:u}:}
First, we introduce a few definitions. Define the \emph{projection-less} derivative at time $t$ for load $i$, $\dot{p}_i(t)$, such that $\dot{x}_i(t)=\Gamma_{\Omega_i,x_i(t)}[\dot{p}_i(t)]$. That is,
\begin{align}\label{eq:pode}
\dot{p}_i(t) \triangleq c\sum\limits_{j\in\N_i}\big(\nabla f_j(x_j(t))-\nabla f_i(x_i(t))\big) + u(t).
\end{align}
Let 
\begin{align}\label{eq:xjstar}
  X_i^*:=\{m \in \R \; |\; \exists \boldx \in X^* \text{ s. t. } \boldx_i = m\}.
\end{align}
\begin{figure}[ht]
  \centering
  \includegraphics[width=0.9\columnwidth]{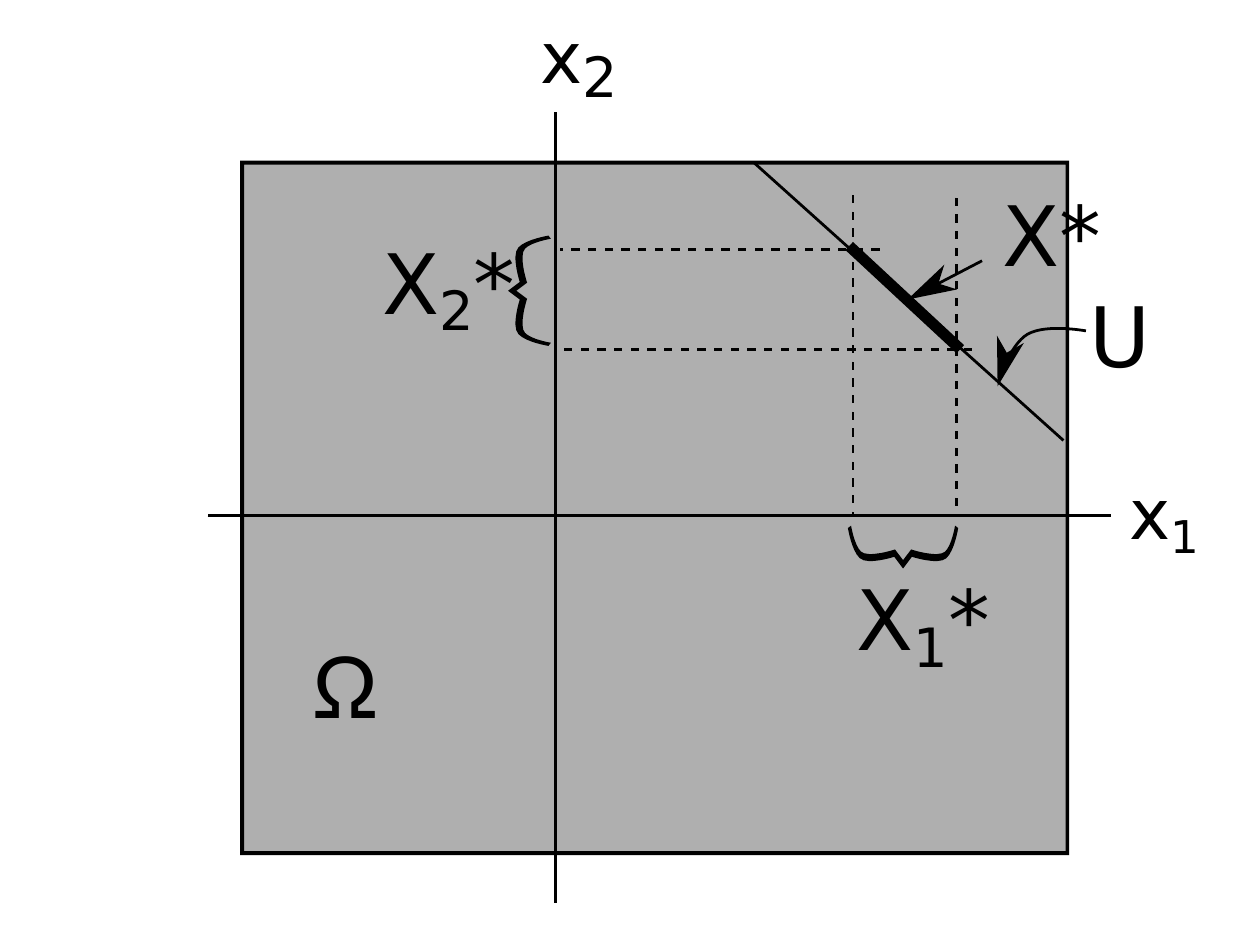}
  \caption{An illustration of the relevant sets for a 2-D case. The equality constraint is satisfied on the line segment, which is the set $U$. The thick sub-segment is the optimal solution set, $X^*$. }
  \label{fig:xjstar-etc}
\end{figure}
Figure \ref{fig:xjstar-etc} shows an illustration for a 2-D case. Note that $x_i\in X_i^*$ for all $i$ does not imply $\boldx\in X^*$ because it may not satisfy the equality constraint, $\boldone^T\boldx=\bar{g}$. However, the sets, $X_i^*$, satisfy a property that is useful in the proofs: if $\boldx \in X_1^* \times X_2^* \times \dots X_n^*$, then $\nabla f_i(x_i) = \nabla f^*$ for every $i$; where $\nabla f^*$ is the optimal gradient; see \eqref{eq:gradfstar}. It should be noted, however, that the converse is true only if $f_i$ is strictly convex. For this reason, we shall call $X_i^*$ the \emph{critical gradient set} of load $i$. Let us define the collections of loads,
\begin{align}\label{eq:ABM}
\begin{aligned}
\A(t)&:=\{i\in\V|x_i(t)>X_i^*\},\\ \B(t)&:=\{i\in\V|x_i(t)<X_i^*\},\\
\M(t)&:=\{i\in\V|x_i(t)\in X_i^*\},
\end{aligned}
\end{align}
At $t$, $\A(t)$ is the set of loads above their critical gradient sets, $\B(t)$ is the set of loads below their critical gradient sets, and $\M(t)$ is the set of loads at their critical gradient sets. By definition, $\A(t)$, $\B(t)$, and $\M(t)$ are mutually disjoint. Define also $a(t):=|\A(t)|$ and $b(t):=|\B(t)|$; note that $a(t)+b(t)\leq n$ for all $t$. Additionally, let
\begin{align}\label{eq:stuck}
\begin{aligned}
\St(t)&:=\{i\in\V|\dot{p}_i(t)\neq\dot{x}_i(t)\},\\
\MM(t)&:=\{i\in\A|\dot{p}_i(t)>0=\dot{x}_i(t)\}\subset\St(t),\\
\Mm(t)&:=\{i\in\B|\dot{p}_i(t)<0=\dot{x}_i(t)\}\subset\St(t).
\end{aligned}
\end{align}
That is, $\St(t)$ is the set of loads for which the projection operation is currently active, $\MM(t)$ is the set of loads in $\A(t)$ currently constrained at their upper bounds due to projection, and $\Mm(t)$ is the set of loads in $\B(t)$ constrained at their lower bounds due to projection. Recall that $\Omega^o$ is the interior of $\Omega$, and note that $X_i^*\subset\Omega^o_i$ because $X^*\subset\Omega^o$; therefore, projection does not affect loads in $\M(t)$. It follows that $\St(t)=\MM(t)\cup\Mm(t)$. Likewise, we have $\MM(t)\subset\A(t)$ and $\Mm(t)\subset\B(t)$. Now, we may observe
\begin{align}
\label{eq:udot}
\begin{aligned}
\dot{u}(t)&=-\boldone^T\dot{\boldx}(t)\\
          &=-\boldone^T\dot{\boldp}(t) + \sum\limits_{i\in\St(t)}\dot{p}_i(t)\\
          &=-nu(t) + \sum\limits_{i\in\MM(t)}\dot{p}_i(t)+\sum\limits_{i\in\Mm(t)}\dot{p}_i(t),
\end{aligned}
\end{align}
where the final equality follows from $\St(t)=\MM(t)\cup\Mm(t)$.

Let the function, $y_i: \R \to \R_{\geq0}$, $i=1,\dots,n$, and the function, $y: \R^n \to \R_{\geq0}$, be defined as
\begin{equation}\label{eq:y}
  \begin{split}
    y_i(m) & \triangleq \inf_r\{|m-r|:r\in X_i^*\}, \quad i\in\V, \\
    y(\boldx) & \triangleq \sum_{i\in\V}y_i(x_i).
  \end{split}
\end{equation}
At time $t$, $y_i(x_i(t))$ is the Euclidean distance between the component, $x_i(t)$, and the critical gradient set, $X_i^*$, and $y(t)$ is the total distance of all loads from their respective critical gradient sets. For the sake of compactness, we will write $y_i(t)$ and $y(t)$ in place of $y_i(x_i(t))$ and $y(\boldx(t))$ in the sequel. Next, we more closely examine the term, $-L\nabla f(\boldx(t))$, in~\eqref{eq:ode}. Without loss of generality, let $\nabla f_1(x_1(t))\leq\nabla f_2(x_2(t))\leq\ldots\leq\nabla f_n(x_n(t))$. 
\begin{align}\label{eq:x_b_increase}
\begin{aligned}
C_B(t)&\triangleq~ c\sum\limits_{i\in\B(t)}\sum\limits_{j\in\N_i}\big(\nabla f_j(x_j(t))-\nabla f_i(x_i(t))\big)\\
&=~c\sum\limits_{i=1}^{b(t)}\sum\limits_{j\in\N_i}\big(\nabla f_j(x_j(t))-\nabla f_i(x_i(t))\big).
\end{aligned}
\end{align}

For $i=1$, every term in~\eqref{eq:x_b_increase} is nonnegative. Suppose $1\notin\N_2$; then for $i=2$, every term in~\eqref{eq:x_b_increase} is nonnegative. Now suppose $1\in\N_2$; then $2\in\N_1$, and the respective terms in~\eqref{eq:x_b_increase} cancel for $i=1,~j=2$ and $i=2,~j=1$. Therefore, the terms of~\eqref{eq:x_b_increase} for $i=1,2$ are nonnegative regardless of whether $1$ and $2$ are neighbors. This argument may be continued for all $i\in\B(t)$. It follows that $C_B(t)\geq0$.

We may define a similar sum over $\A(t)$:
\begin{align}\label{eq:x_a_decrease}
C_A(t)\triangleq&~ c\sum\limits_{i\in\A(t)}\sum\limits_{j\in\N_i}\big(\nabla f_j(x_j(t))-\nabla f_i(x_i(t))\big).
\end{align}
By the same argument as above, we may show that $C_A(t)\leq0$.

Combining~\eqref{eq:pode} and~\eqref{eq:x_b_increase}, we can see that the total change in consumption for loads in $\B(t)$ \emph{without projection} is
\begin{align*}
\sum\limits_{i\in\B(t)}\dot{p}_i(t)=C_B(t)+b(t)u(t).
\end{align*}
Similarly from~\eqref{eq:pode} and~\eqref{eq:x_a_decrease}, the total change in consumption for loads in $\A(t)$ without projection is
\begin{align*}
\sum\limits_{i\in\A(t)}\dot{p}_i(t)=C_A(t)+a(t)u(t).
\end{align*}
%and the total change in consumption for loads in $\M(t)$ without projection is
%\begin{align*}
%\sum\limits_{i\in\M(t)}\dot{p}_i(t)=-C_A(t)-C_B(t)+m(t)u(t).
%\end{align*}

Note that
\begin{align*}
\sum\limits_{i\in\B(t)}\dot{x}(t)&=\sum\limits_{i\in\B(t)}\big(\dot{p}_i(t)\big) - \sum\limits_{i\in\Mm(t)}\big(\dot{p}_i(t)\big)\\
&=C_B(t) + b(t)u(t) - \sum\limits_{i\in\Mm(t)}\dot{p}_i(t)
\end{align*}
and
\begin{align*}
\sum\limits_{i\in\A(t)}\dot{x}(t)&=\sum\limits_{i\in\A(t)}\big(\dot{p}_i(t)\big) - \sum\limits_{i\in\MM(t)}\big(\dot{p}_i(t)\big)\\
&=C_A(t) + a(t)u(t) - \sum\limits_{i\in\MM(t)}\dot{p}_i(t).
\end{align*}
Because $\boldx(t)$ is a Caratheodory solution (see Proposition~\ref{prop:nag}), it is absolutely continuous and differentiable almost everywhere. It follows that $y(t)$ is absolutely continuous and differentiable almost everywhere. Then we have
\begin{align}\label{eq:ydotlong}
\begin{aligned}
\dot{y}(t)&=\sum\limits_{i\in\A(t)}\big(\dot{x}_i(t)\big) - \sum\limits_{i\in\B(t)}\big(\dot{x}_i(t)\big)\\
          &=\Big(C_A(t)+a(t)u(t)-\sum\limits_{i\in\MM(t)}\dot{p}_i(t)\Big)\\
          &\quad - \Big(C_B(t)+b(t)u(t)-\sum\limits_{i\in\Mm(t)}\dot{p}_i(t)\Big)\\
          &=C_A(t)-C_B(t) + \big(a(t)-b(t)\big)u(t)\\
          &\quad - \sum\limits_{i\in\MM(t)}\big(\dot{p}_i(t)\big)+ \sum\limits_{i\in\Mm(t)}\big(\dot{p}_i(t)\big).
\end{aligned}
\end{align}
Observe that, for $u(t)\geq0$, we have
\begin{align}\label{eq:y+u}
\dot{y}(t) + \dot{u}(t) &= C_A(t)-C_B(t) +\big(a(t)-b(t)-n\big)u(t) \notag \\
                        &\quad \quad +2\sum_{i\in\Mm(t)} \dot{p}_i(t)\\
                        &\leq -u(t) = -|u(t)| \leq0, \notag
\end{align}
where we have used $a(t)\leq n-1$ if $u(t)\geq0$ (i.e., if $u(t)\geq0$, at least one load $i$ must be at or below it critical gradient set, $X_i^*$, by definition of $u(t)$). 
Similarly, for $u(t)\leq0$, we have
\begin{align}\label{eq:y-u}
\dot{y}(t) - \dot{u}(t) &= C_A(t)-C_B(t) + \big(a(t)-b(t)+n\big)u(t) \notag\\
                        &\quad - 2\sum\limits_{i\in\MM(t)}\dot{p}_i(t)\\
                        &\leq u(t) =-|u(t)| \leq0. \notag
\end{align}
Consider the function
\begin{align}\label{eq:z}
  z: \R^n \to \R_{\geq0}: z(\boldx) & = y(\boldx)+|\bar{g} - \boldone^T\boldx| \\
 \Rightarrow   z(t) & = y(t)+|u(t)|. \notag
\end{align}
Just as $y(t)$ is absolutely continuous and differentiable almost everywhere, so is $|u(t)|$, and therefore, $z(t)$. Then, we may combine~\eqref{eq:y+u} and~\eqref{eq:y-u} to obtain
\begin{align}\label{eq:zdot}
\dot{z}(t) \leq -|u(t)| \leq 0.
\end{align}
%\begin{align*}
%\dot{z}(t) = \left\{\begin{array}{l}
%\dot{y}(t) + \dot{u}(t), \quad u(t)>0\\
%\dot{y}(t) - \dot{u}(t), \quad u(t)<0\\
%\dot{y}(t) + |\dot{u}(t)|, \;\; u(t)=0. \end{array}\right.
%\end{align*}
%It follows from~\eqref{eq:ydotlong},~\eqref{eq:y+u}, and~\eqref{eq:y-u} that 
Hence, $z(t)$ is non-increasing almost everywhere.

Now, suppose $u(t)$ does not converge to 0. For some $t_1<t_2$, let $|u(t)|>\mu$ for $t_1<t<t_2$ and some $\mu>0$. From~\eqref{eq:zdot} and the first fundamental theorem of calculus, it then follows that
\begin{align*}
z(t_2)-z(t_1)\leq-\int\limits_{t_1}^{t_2}\mu~dt.
\end{align*}
Because $u(t)$ is continuous and $u(t)\nrightarrow 0$, there exists a sequence, $t_1,t_2,\ldots,$ with $t_k \to \infty$ as $k \to \infty$ such that $|u(t)|>\mu$ for all $t_1<t<t_2,~t_3<t<t_4,~\ldots$. During all other intervals (e.g., $[t_2,t_3]$), %At all other $t$'s in $[t_1, ~t_k]$, 
$z(t)$ does not increase since $\dot{z}(t)\leq 0$ for all $t$. This implies $z(t)$ decreases without bound, which is impossible since the boundedness of the domain, $\Omega$, implies $z(t)$ must be bounded. Therefore, $u(t)\to0$.\hfill$\square$\\

\noindent\textbf{Proof of Lemma \ref{lem:stable}:} 
Global asymptotic stability of an invariant set requires stability and global attractiveness (see Definition 8.1 in section 8.4 in~\cite{Khalil}). Theorem~\ref{the:main} establishes global attractiveness of $X^*$, so only stability remains to be proven. For stability, it is sufficient to show that each neighborhood of $X^*$ is positively invariant~\cite{Khalil}. In the proof of Lemma~\ref{lem:u}, it was shown that $z(t)$ is non-increasing~\eqref{eq:zdot}, where $z(t)$ is defined in~\eqref{eq:z}. Let $Z_m\triangleq\{\boldx\in\Omega|z<m\}$. It follows from the definition that any neighborhood of $X^*$ takes the form of $\mathcal{Z}_m$ for some $m$~\cite{Khalil}. Since $z(t)$ is non-increasing, it follows that $\mathcal{Z}_m$ is a positively invariant neighborhood of $X^*$ for each $m\geq0$. %Let $\mathcal{O}$ be any neighborhood of $X^*$, and let $Z:=\min_{\boldx}\{z|\boldx\in\mathcal{O}\}$. Then $\mathcal{Z}_Z$ is a positively invariant neighborhood of $X^*$ contained in the neighborhood, $\mathcal{O}$. 
Therefore, $X^*$ is stable. This concludes the proof.\hfill$\square$

\section{Simulation results}\label{sec:results}
\subsection{Simulation setup}\label{sec:sim}
\begin{figure}[htpb]
\begin{center}
\includegraphics[scale=0.25]{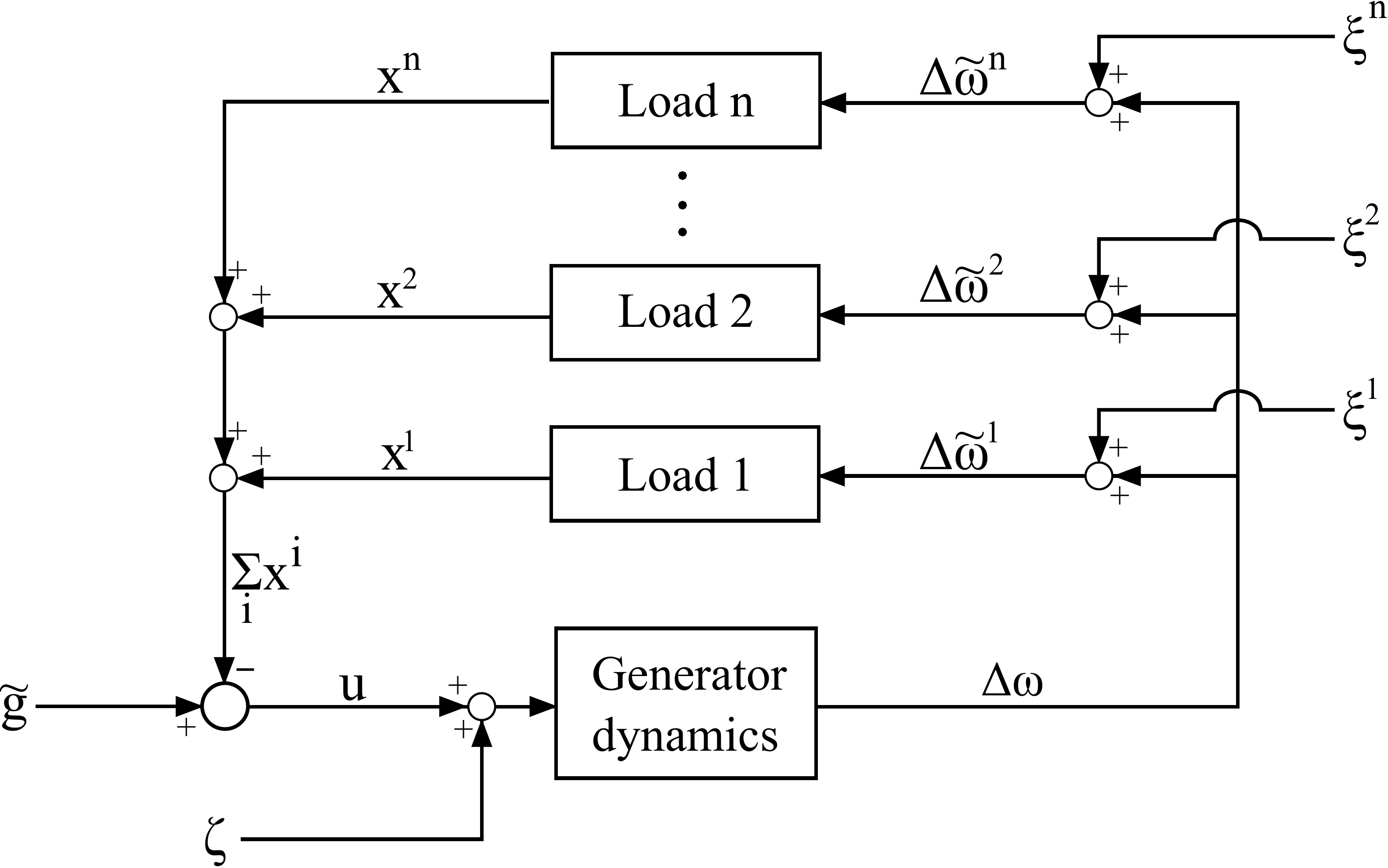}
\caption{System architecture for simulations. Inter-load communication is not shown.}
\label{fig:powergrid}
\end{center}
\end{figure}

Figure~\ref{fig:powergrid} shows the system architecture used for design and simulation. 
The generator dynamics block shown in Figure~\ref{fig:powergrid} also includes local controls that are usually present in generators. The loss of generation is modeled as an exogenous disturbance, $\bar{g}$, in the figure. The estimator in the figure is the one described in Section~\ref{sec:estimator} to estimate the consumption-generation mismatch from local, noisy frequency measurements. The process disturbance, $\zeta$, and measurement noise, $\xi_i$, at each load are modeled as wide-sense stationary white noise. %, and $\xi^i$ is independent across loads. 

For ease of comparison between the DGP algorithm and the dual algorithm, we use the same generator dynamics, noise statistics, and communication graph as in~\cite{zhao2013optimal}, and the reader is referred to that work for more detailed information about the simulation model or implementation of the state estimator.

Even without the use of intelligent loads, the local generator control will change the generator set point in response to frequency deviation to match consumption, which will restore the frequency to its nominal value on its own. Intelligent loads are supposed to help the generator in reacting to frequency deviations faster so that large excursions of system frequency are avoided.

%To estimate $u$, loads use a state observer described in~\cite{kitanidis1987unbiased}, where the frequency deviation, $\Delta\omega$, is the state. Each load $i$ uses the previously predicted state, $\Delta\hat{\omega}^i_{t|t-1}$, and the measured $\Delta\tilde{\omega}^i_t$ to obtain $\hat{u}_t$, and the estimation error forms a martingale-difference 
% sequence~\cite{zhao2013optimal}.

For each load $i$, we choose $\Omega_i=[-\bar{x}_i,\bar{x}_i]$, where $\bar{x}_i$ is chosen from a uniform distribution and then normalized so that $\sum_{i=1}^{n}\bar{x}_i=60$ MW (as in~\cite{zhao2013optimal}).

We test the performance of the DGP algorithm with two distinct disutility functions. The first is a convex but not strictly convex function:
\begin{align}\label{eq:f-case2}
  f_i(x_i) =
  \begin{cases}
    0, & |x_i| < a_i\\
    q_i(x_i-a_i)^2, & x_i \geq a_i \\ 
    q_i(x_i+a_i)^2, & x_i \leq -a_i
\end{cases}
\end{align}
where $a_i=0.1\bar{x}_i$. The consumer does not experience any disutility as long as the load variation is within $\pm a_i$. The second disutility function \emph{is} strictly convex: 
\begin{align}\label{eq:quad-f}
f_i(x_i) = q_i(x_i)^2.
\end{align}
For both disutility functions, we pick $q_i$ to be an arbitrary positive number such that $1/q_i$ is chosen from a uniform distribution on the interval of $[0.1,\;0.3]$. This is chosen for comparison with~\cite{zhao2013optimal}, which makes a similar choice for disutility functions.

The initial conditions are $g[0] = 200$ MW and $u[0] = 0$. Two generation contingencies are modeled as step changes:
\begin{equation*}
g[k] = \left\{\begin{array}{l}
200~\mathrm{MW},\;\;~0~\mathrm{s} \leq kT < 20~\mathrm{s}\\
190~\mathrm{MW},~20~\mathrm{s} \leq kT < 50~\mathrm{s}\\
170~\mathrm{MW},~50~\mathrm{s} \leq kT, \end{array} \right.
\end{equation*}
where $T=0.1$ seconds is the discretization interval.

\begin{figure}[ht]
  \centering
  \includegraphics[width=0.9\columnwidth]{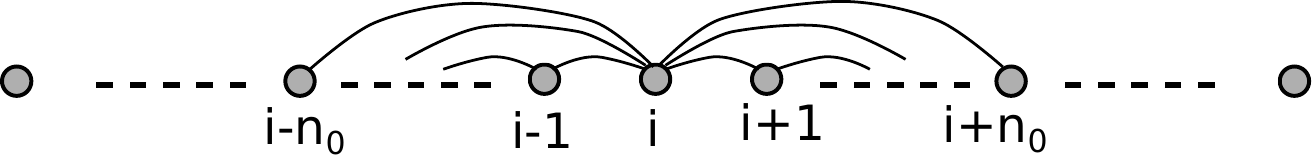}
  \caption{The communication graph. To avoid clutter, not all edges are shown.}
  \label{fig:network}
\end{figure}

%As in~\cite{zhao2013optimal}, t
Simulations are conducted with the communication network in~\cite{zhao2013optimal}, where load $i$ communicates with loads from max\{1, $i-n_0$\} to min\{$n$, $i+n_0$\}, where $n_0\leq n$. We use $n=1000$ and $n_0=1$. The network is shown in Figure \ref{fig:network}.

Additionally, we use $c=5$ and $\gamma[k]=\gamma[0]/(k^{0.8})$ for $k>0$, with $\gamma[0]=1.5\underline{q}/n$, where $\underline{q}\triangleq\min_i{q_i}$.

With all of these parameter choices, Assumptions~\ref{as:f},~\ref{as:power}, and~\ref{as:tech} are satisfied. Note that Assumption~\ref{as:tech}(3) is satisfied from the discussion in Section~\ref{sec:estimator}. %Theorem~\ref{the:main} is applicable to both disutility functions~\eqref{eq:f-case2} and~\eqref{eq:quad-f}. Theorem~\ref{the:quadratic} applies to disutility function~\eqref{eq:quad-f}.

\subsection{Results with non-strictly convex disutility}
\begin{figure}[htpb]
\begin{center}
\includegraphics[scale=0.3]{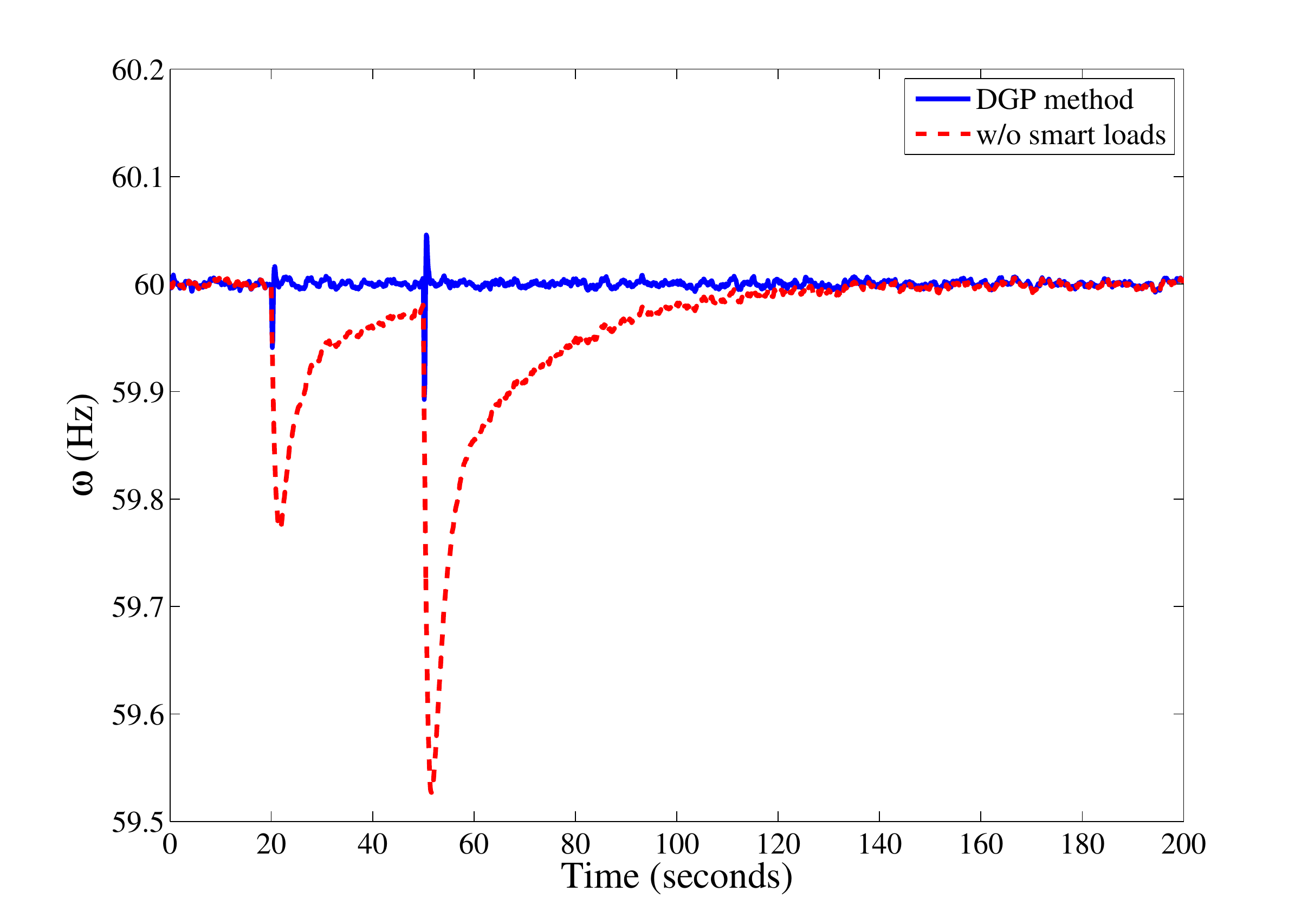}
\includegraphics[scale=0.3]{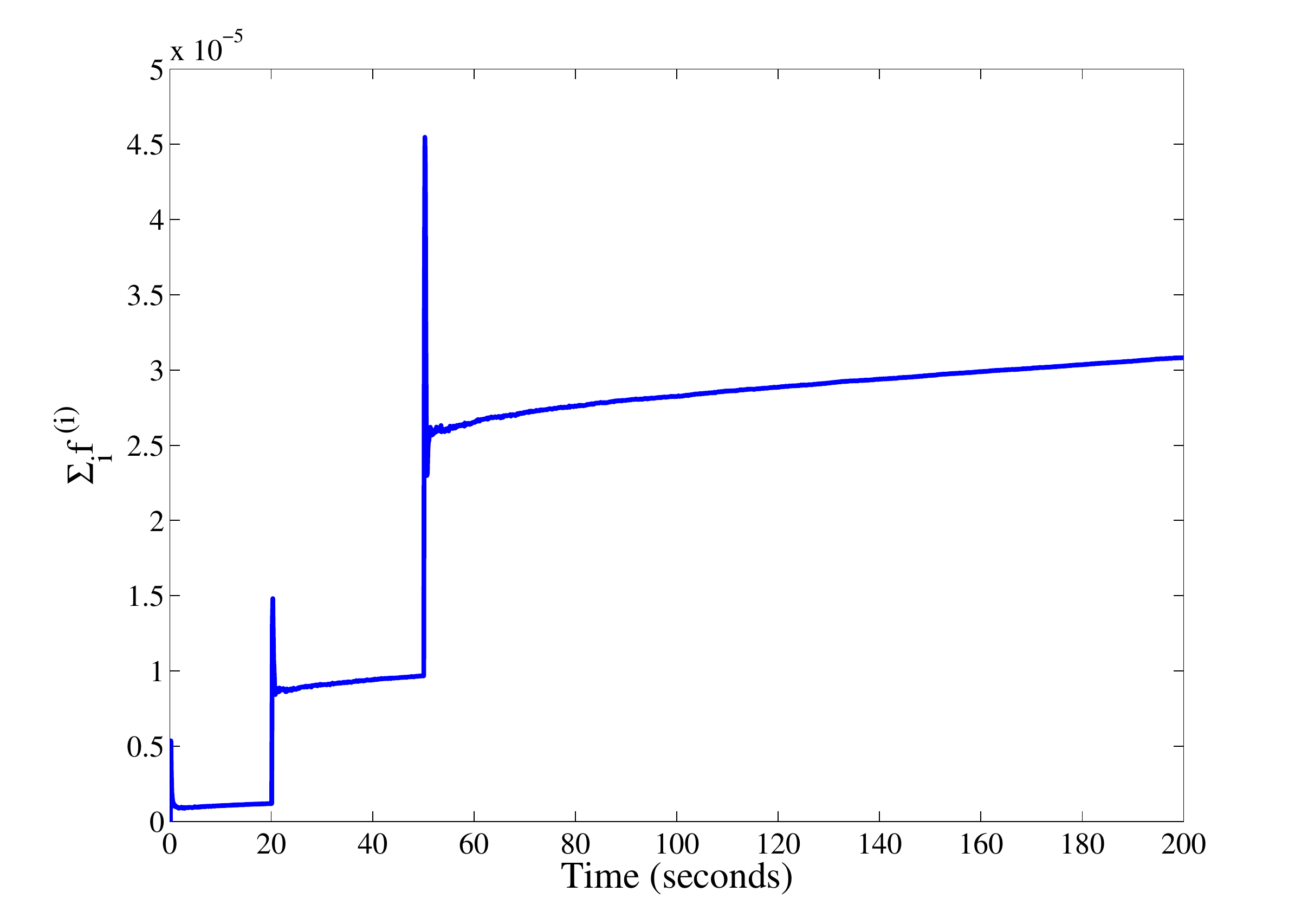}
\caption{Performance of the DGP algorithm with consumer disutility that is not strictly convex. Step changes in generation occur at 20 and 50 seconds.}
\label{fig:flatf_nopro_1000loads}
\end{center}
\end{figure}

Figure~\ref{fig:flatf_nopro_1000loads} shows simulation results for the DGP algorithm with consumer disutility function~\eqref{eq:f-case2}. The dual algorithm (DA) from~\cite{zhao2013optimal} is not applicable because the inverse of $\nabla f(\boldx)$ must exist in $\Omega$ to implement DA, which is not the case when $|x_i| \leq a_i$.

The system frequency without smart loads (i.e., with generator-only control) is shown in red as well. 
Using DGP, the loads are able to assist the generator in avoiding large frequency deviations from the nominal when each contingency occurs.

\subsection{Comparison with dual algorithm: results with strictly convex disutility}

\begin{figure}[htpb]
\begin{center}
\includegraphics[scale=0.33]{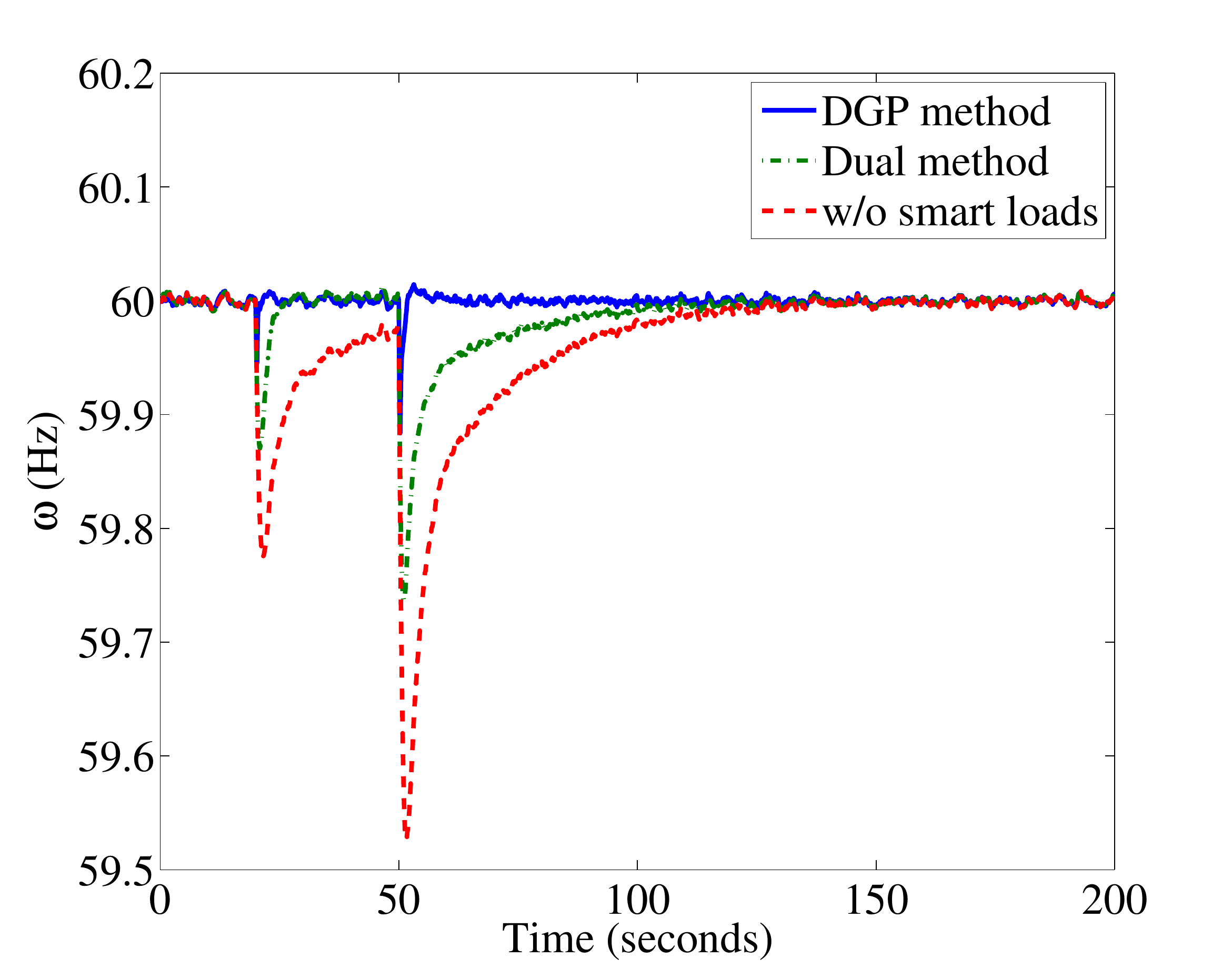}
\includegraphics[scale=0.3]{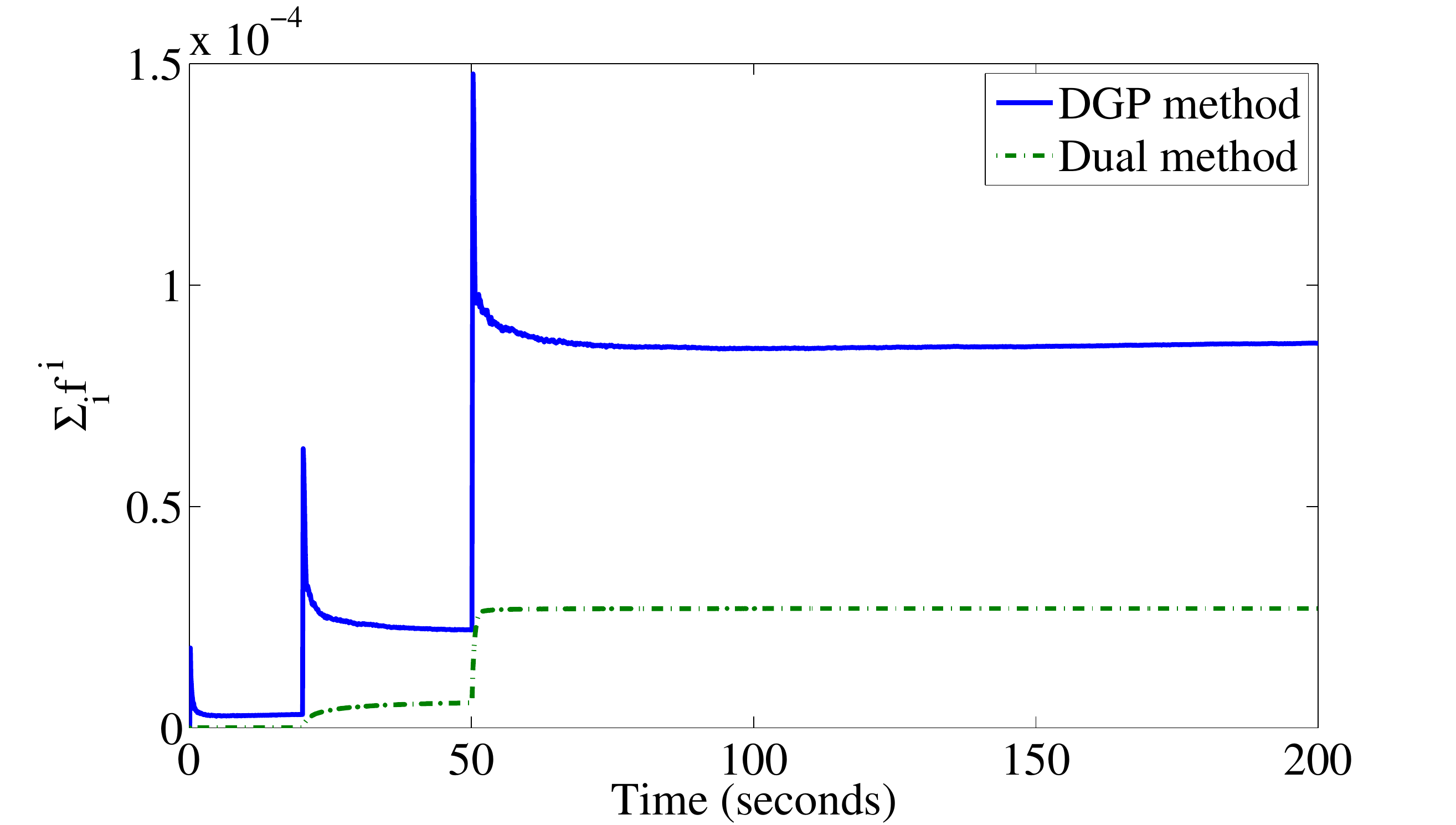}
\caption{Performance of the DGP and dual algorithms with quadratic consumer disutility with projection. Step changes in generation occur at 20 and 50 seconds.}
\label{fig:sync_freq_1000loads}
\end{center}
\end{figure}

Figure~\ref{fig:sync_freq_1000loads} shows results of DGP and DA with quadratic disutilities~\eqref{eq:quad-f}.
%In both the DGP and dual algorithms, the diminishing $\gamma_k$ may slow down response. If the time of each generation drop were known exactly, $\gamma_k$ may be reinitialized at those times---resulting in faster response. 
% when $\gamma_k$ is reinitialized at the time of each generation drop. That is, at $r=20$ seconds and $s=50$ seconds we set $\gamma_r=\gamma_{k-r}$ for $r\leq k < s$ and $\gamma_s=\gamma_{k-s}$ for $s\leq k$. This can be viewed as the ``ideal'' scenario.
DGP results in a significantly smaller frequency drop compared to both generator-only control and DA.
Although DA returns the frequency to the nominal value faster than generator-only control, it does not reduce the initial frequency drop as much as DGP.

However, the consumer disutility is significantly lower for DA than for DGP. This is because DA is responding more slowly than DGP, so the equality constraint is not being satisfied---resulting in a lower cost. The slower response of DA is due to the inversion of the derivative of each load's disutility function. Because the derivative of each disutility function is rather steep, the inverse is quite flat, so large changes in its argument still result in small changes in its value---leading to small changes in consumption. Conversely, DGP aggressively meets the equality constraint because of the generation-matching step. This results in a lower frequency deviation but more disutility.

%DA has a significantly lower steady-state disutility because the generator control restores much of the frequency. The loads interpret the restored frequency as a smaller consumption-generation mismatch, which results in less change in consumption and therefore lower disutility.

{\bf Effect of network size and structure:} Although not reported here, simulations with varying number of loads ($n=10,~100$) and varying amount of communication ($n_0=10,~100,~1000$) showed similar trends as in the $n=1000$, $n_0=1$ case. 
It was observed in~\cite{zhao2013optimal} that DA showed similar behavior. The reason for this insensitivity to network size and structure is likely the use of local frequency measurements which provides global information, effectively creating a virtual communication link between a fictitious central node and all nodes of the network.

\section{Conclusion}
\label{sec:conc}
The DGP algorithm solves a constrained optimization problem in a distributed manner to aid a power grid in maintaining system frequency near its nominal value while minimizing consumers' disutility. The main advantage of the DGP algorithm over the dual algorithm of~\cite{zhao2013optimal} is that it is applicable to disutility functions that are convex but not necessarily strictly convex, while the dual algorithm can be used only for strictly convex functions. Disutility functions that are not strictly convex model more realistic consumer behavior that is insensitive to small changes in consumption. Simulations also show that the DGP algorithm performed either better than or similarly to the dual algorithm from~\cite{zhao2013optimal} in maintaining frequency. 

The convergence of the DGP algorithm required that the optimal points lie in the strict interior of the domain.  When capacity of the loads is small, the optimal solution is likely to lie on the boundary. An open problem is the design of a distributed primal algorithm, if such an algorithm exists, that is guaranteed to converge to an optimal solution lying on the boundary.

The analysis in this paper assumed a time-invariant communication graph. In contrast, the dual algorithm in~\cite{zhao2013optimal} was proved to converge even with a time-varying communication graph. Convergence analysis of the DGP algorithm for the time-varying case is left for future work.

% The work reported here assumes the power grid is a microgrid, since the grid-frequency is assumed constant throughout the grid. Future work will focus on extending the DGP algorithm into more-realistic power grid models in which the frequency varies both spatially and temporally. Other interesting avenues fo future work include extension of the DGP algorithm to time-varying communication networks. 

One issue that was ignored in our analysis---as well as that in~\cite{zhao2013optimal}---is the effect of feedback interconnection between the generator control system and the load control algorithm: % DGP in this paper and the dual algorithm in~\cite{zhao2013optimal}. 
simulations reported here and in~\cite{zhao2013optimal} are conducted with both the control systems in place, and simulation results do not show instability. However, analysis to rule out possible instabilities is lacking. A related issue is actuator dynamics. It is assumed that loads can react as fast as asked by the load control algorithm. The phase lag due to loads' inertia can reduce performance or even cause instability. These are interesting topics worthy of future investigation.

\section*{Acknowledgment}
The authors thank C. Zhao and S. Low for their assistance in reproducing the results of~\cite{zhao2013optimal}, Sean Meyn for pointing out the connection to the Skorokhod problem, and Vivek Borkar for several useful comments, including pointing us to the work by Nagurney. The  research reported here was partially supported by the National Science Foundation through grants 1463316 and 1646229.

\bibliographystyle{plain}

% JB
\bibliography{../brooks_cloud/bibs/Jonathan,../PBbib/building,../PBbib/peoplemodel,../PBbib/Barooah,../PBbib/sensnet_bib_dbase,../PBbib/optimization,../PBbib/grid}

\begin{thebibliography}{10}

\bibitem{bolzam:13}
S.~Bolognani and S.~Zampieri.
\newblock A distributed control strategy for reactive power compensation in
  smart microgrids.
\newblock {\em {IEEE} Trans. on Automatic Control}, 58(11), November 2013.

\bibitem{borkarbook:2008}
Vivek~S. Borkar.
\newblock {\em Stochastic Approximation: A Dynamical Systems Viewpoint}.
\newblock Cambridge University Press, 2008.

\bibitem{brobar:ACC:2016}
Jonathan Brooks and Prabir Barooah.
\newblock Consumer-aware load control to provide contingency reserves using
  frequency measurements and inter-load communication.
\newblock In {\em {A}merican Control Conference}, pages 5008 -- 5013, July
  2016.

\bibitem{cherukuri2014initialization}
Ashish Cherukuri and Jorge Cortes.
\newblock Initialization-free distributed coordination for economic dispatch
  under varying loads and generator commitment.
\newblock {\em Automatica}, 2014.
\newblock submitted.

\bibitem{cor:2008}
Jorge Cort\'{e}s.
\newblock Discontinuous dynamical systems - a tutorial on solutions, nonsmooth
  analysis, and stability.
\newblock {\em {IEEE Control Systems Magazine}}, 28(3):36--73, 2008.

\bibitem{GodsilRoyle_2001}
Chris Godsil and Gordon Royle.
\newblock {\em Algebraic Graph Theory}.
\newblock Graduate Texts in Mathematics. Springer, 2001.

\bibitem{khadgi2014modeling}
Prajwal Khadgi, Lihui Bai, and Gerald Evans.
\newblock Modeling demand response using utility theory and model predictive
  control.
\newblock In {\em IIE Annual Conference. Proceedings}, page 1262. Institute of
  Industrial Engineers-Publisher, 2014.

\bibitem{Khalil}
H.K. Khalil.
\newblock {\em {Nonlinear Systems 3rd}}.
\newblock Prentice hall Englewood Cliffs, NJ, 2002.

\bibitem{kirby2007ancillary}
Brendan Kirby.
\newblock Ancillary services: Technical and commercial insights.
\newblock 2007.
\newblock prepared for W\"{a}rtsil\"{a} North America Inc.

\bibitem{kitanidis1987unbiased}
Peter~K Kitanidis.
\newblock Unbiased minimum-variance linear state estimation.
\newblock {\em Automatica}, 23(6):775--778, 1987.

\bibitem{kushner2003stochastic}
Harold Kushner and G~George Yin.
\newblock {\em Stochastic approximation and recursive algorithms and
  applications}, volume~35.
\newblock 2nd edition, 2003.

\bibitem{lew2013western}
Debra Lew, Greg Brinkman, E~Ibanez, BM~Hodge, and J~King.
\newblock The western wind and solar integration study phase 2.
\newblock {\em National Renewable Energy Laboratory, NREL/TP-5500}, 55588,
  2013.

\bibitem{linbarmeymid:2015}
Yashen Lin, Prabir Barooah, Sean Meyn, and Timothy Middelkoop.
\newblock Experimental evaluation of frequency regulation from commercial
  building {HVAC} systems.
\newblock {\em {IEEE} Transactions on Smart Grid}, 6:776 -- 783, 2015.

\bibitem{Luenbergerbluebook:2003}
David~G. Luenberger.
\newblock {\em Linear and Nonlinear Programming}.
\newblock Springer, 2 edition, 2003.

\bibitem{pnnl2008value}
Y.~V. Makarov, Lu. S., J.~Ma, and T.~B. Nguyen.
\newblock Assessing the value of regulation resources based on their time
  response characteristics.
\newblock Technical Report PNNL-17632, Pacific Northwest National Laboratory,
  Richland, WA, June 2008.

\bibitem{milligan2010utilizing}
Michael Milligan and Brendan Kirby.
\newblock Utilizing load response for wind and solar integration and power
  system reliability.
\newblock In {\em Wind Power Conference, Dallas, Texas}, 2010.

\bibitem{molina2011decentralized}
Angel Molina-Garc{\'\i}a, Fran{\c{c}}ois Bouffard, and Daniel~S Kirschen.
\newblock Decentralized demand-side contribution to primary frequency control.
\newblock {\em Power Systems, IEEE Transactions on}, 26(1):411--419, 2011.

\bibitem{nagurney2012projected}
Anna Nagurney and Ding Zhang.
\newblock {\em Projected dynamical systems and variational inequalities with
  applications}, volume~2.
\newblock Springer Science \& Business Media, 2012.

\bibitem{nedic2009approximate}
Angelia Nedic and Asuman Ozdaglar.
\newblock Approximate primal solutions and rate analysis for dual subgradient
  methods.
\newblock {\em SIAM Journal on Optimization}, 19(4):1757--1780, 2009.

\bibitem{NERCBalancingDocument:2011short}
{Prepared by NERC RS Committee}.
\newblock {Balancing and Frequency Control: A Technical Document Prepared by
  the NERC Resources Subcommittee}.
\newblock {NERC Technical Report}, January 26 2011.

\bibitem{schFAPER80}
F.C. Schweppe, R.D. Tabors, J.L. Kirtley, H.R. Outhred, F.H. Pickel, and A.J.
  Cox.
\newblock Homeostatic utility control.
\newblock PAS-99(3):1151 --1163, May 1980.

\bibitem{short2007stabilization}
Joe~A Short, David~G Infield, and Leon~L Freris.
\newblock Stabilization of grid frequency through dynamic demand control.
\newblock {\em Power Systems, IEEE Transactions on}, 22(3):1284--1293, 2007.

\bibitem{siano2014demand}
Pierluigi Siano.
\newblock Demand response and smart grids: a survey.
\newblock {\em Renewable and Sustainable Energy Reviews}, 30:461--478, 2014.

\bibitem{skorokhod1961stochastic}
Anatoliy~V Skorokhod.
\newblock Stochastic equations for diffusion processes in a bounded region.
\newblock {\em Theory of Probability \& Its Applications}, 6(3):264--274, 1961.

\bibitem{todd2008providing}
D~Todd, M~Caufield, B~Helms, Alcoa~Power Generating, Inc~M Starke, B~Kirby, and
  J~Kueck.
\newblock Providing reliability services through demand response: A preliminary
  evaluation of the demand response capabilities of {ALCOA} {INC.}
\newblock {\em ORNL/TM}, 233, 2008.

\bibitem{zhang2014optimal}
Bowen Zhang, Michael~C Caramanis, and John Baillieul.
\newblock Optimal price-controlled demand response with explicit modeling of
  consumer preference dynamics.
\newblock In {\em Decision and Control (CDC), 2014 IEEE 53rd Annual Conference
  on}, pages 2481--2486. IEEE, 2014.

\bibitem{zhao2013optimal}
Changhong Zhao, Ufuk Topcu, and Steven~H Low.
\newblock Optimal load control via frequency measurement and neighborhood area
  communication.
\newblock {\em Power Systems, IEEE Transactions on}, pages 3576--3587, 2013.

\bibitem{zhu2012distributed}
Minghui Zhu and Sonia Mart{\'\i}nez.
\newblock On distributed convex optimization under inequality and equality
  constraints.
\newblock {\em Automatic Control, IEEE Transactions on}, 57(1):151--164, 2012.

\bibitem{zhu2013approximate}
Minghui Zhu and Sonia Mart{\'\i}nez.
\newblock An approximate dual subgradient algorithm for multi-agent non-convex
  optimization.
\newblock {\em Automatic Control, IEEE Transactions on}, 58(6):1534--1539,
  2013.

\end{thebibliography}
%PB
%\bibliography{IEEEabrv,../../bibs/Jonathan,../../PBbib/building,../../PBbib/sensnet_bib_dbase,../../PBbib/grid,../../PBbib/Barooah,../../PBbib/optimization}
\end{document}